\documentclass[12pt]{article}
\usepackage{amsmath,amssymb,color}
\usepackage{cite,epsf} 
\usepackage{pstricks}
\usepackage{color}
\usepackage{hyperref}
\usepackage{tikz}     

\usetikzlibrary{
decorations.pathmorphing 
}

\usepackage{graphicx,array} 
\newcommand{\be}{\begin{equation}}
\newcommand{\ee}{\end{equation}}
\newcommand{\beq}{\begin{equation}}
\newcommand{\eeq}{\end{equation}}
\newcommand{\bea}{\begin{eqnarray}}
\newcommand{\eea}{\end{eqnarray}}

\newcommand{\ba}{\begin{eqnarray}}
\newcommand{\ea}{\end{eqnarray}}

\def\cos{\mbox{cos}}

\begin{document}

\begin{titlepage}
\vspace{10pt}
\hfill
{\large\bf HU-EP-24/03}
\vspace{20mm}
\begin{center}

{\Large\bf  Conformal invariants of curves  via those\\[2mm]  for inscribed polygons with circular edges
}

\vspace{45pt}

{\large Harald Dorn 
{\footnote{dorn@physik.hu-berlin.de
 }}}
\\[15mm]
{\it\ Institut f\"ur Physik und IRIS Adlershof, 
Humboldt-Universit\"at zu Berlin,}\\
{\it Zum Gro{\ss}en Windkanal 6, D-12489 Berlin, Germany}\\[4mm]

\vspace{20pt}

\end{center}
\vspace{10pt}
\vspace{40pt}

\centerline{{\bf{Abstract}}}
\vspace*{5mm}
\noindent 
The conformal nature of  smooth curves in  $\mathbb{R}^3$ is characterised by conformal length, curvature and torsion. We present a derivation of these conformal parameters via a limiting process using inscribed polygons with circular edges .
The procedure is based on elementary geometry in $\mathbb{R}^3$ only and similar to the rectification of curves in the metrical case. It seems
to be not available in the literature so far.
\end{titlepage}
\section{Introduction}
The conformal geometry of smooth curves in $\mathbb{R}^3$ and in higher dimensions has been studied by mathematicians using various techniques, as e.g. 
conformal derivatives, normal forms, group theoretical methods or the kinematics of osculating spheres, see e.g. \cite{liebmann,Fialkow,Cairns,Sulanke,Fuster1,Fuster2,Langevin} and references therein. \\
In analogy to the metrical invariants length $s$, curvature $\kappa$ and torsion $\tau$ one  finds for smooth 3-dimensional curves,  see e.g. \cite{Cairns},
conformal length $\omega$,\footnote{The prime denotes derivative with respect to $s$.}
\beq
d\omega~=~\sqrt{\nu}ds,~~~~~~\nu~=~\sqrt{(\kappa ')^2+\kappa ^2\tau ^2}~,\label{length-omega}
\eeq
conformal curvature $Q$
\beq
Q~=~\frac{4(\nu ''-\kappa ^2\nu)\nu -5(\nu ')^2}{8\nu ^3}~,\label{Q}
\eeq
and conformal torsion $T$
\beq
T~=~\frac{2(\kappa ')^2\tau+\kappa^2\tau ^3+\kappa\kappa'\tau '-\kappa\kappa''\tau}{\nu ^{\frac{5}{2}}}~.\label{T}
\eeq

Our interest in this subject  originates from  studies of Wilson loops \cite{Dorn:2020meb,Dorn:2020vzj,Dorn:2023qbj} related to polygons with circular edges.
There the conformal length of the edges is zero and both $Q$ and $T$ are ill-defined along the edges. As conformal invariant parameters fixing these polygons,
up to conformal transformations, we identified the cross ratios formed out of the positions of the corners, the angles between the edges (cusp angles) and angles between circles fixed by three consecutive corners and adjacent edges (called torsion angles). For these polygons conformal curvature and torsion is somehow
located at the corners in a distributional sense. So far we could not give this comment a precise meaning.

But instead we were able to realise the
opposite direction, i.e. to relate $\omega,Q$ and $T$ for smooth curves to cusp and torsion angles of inscribed polygons.\footnote{We did not find
in the literature such a direct and vivid derivation by means of elementary geometry  in $\mathbb{R}^3$.  Therefore, we decided to write it down in this note.} The procedure is analogously to the
derivation of metrical length, curvature and torsion via rectification by ordinary polygons with straight edges. Such a rectification has also the potential
to be applied to more general irregular curves, for the metrical case see e.g.\cite{book}. Rectifications by ordinary polygons are not conformal covariant, since in general straight lines are not mapped to straight lines. But circles are mapped to circles, hence one can glue together pieces of successive circles fixed by sets of three nearest neighbour points to a conformal covariant approximation (rectification) of a given curve, see fig.1.
\begin{figure}[h!]
\begin{center}
  \includegraphics[width=11cm]{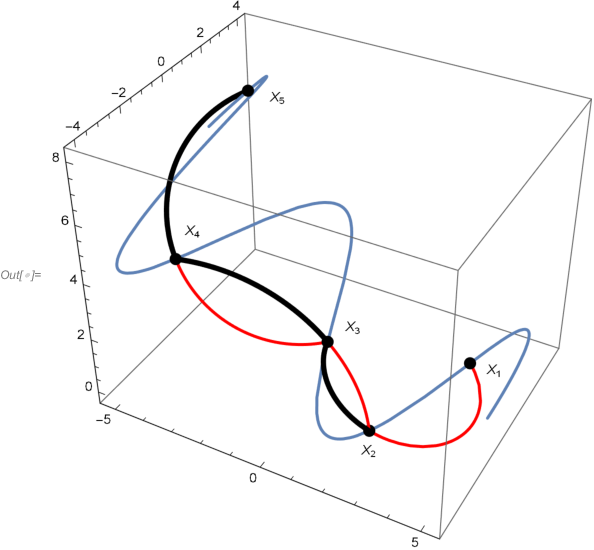}
  \end{center}

\caption {\it In blue is shown a certain smooth curve with five randomly chosen points $x_1,\dots ,x_5$. Each set of three consecutive
points   $(x_{j-1},x_j,x_{j+1})$ fixes a circle $cc_j$. The arcs from $x_{j-1}$ to $x_j$ are depicted in red, those from $x_j$ to $x_{j+1}$ in black. Both the black and red polygon is a conformal covariant approximation (rectification) of the respective piece of the blue smooth curve.}
\label{rectif} 
\end{figure}
The paper is organised as follows. To emphasize the analogy to the metrical case, section 2 gives a short recapitulation of curvature and torsion for smooth curves via a limiting process with inscribed ordinary polygons. The main result for the conformal invariants via inscribed polygons with circular edges is derived
in section 3. Finally in section 4 we add some by the way obtained observations on cross ratios and circumcircles of triangles. An appendix lists formulas
from elementary geometry concerning circumcircles and circumspheres.
\section{A short recap of the metrical case}
Let denote by $x(s)$ a smooth curve in  $\mathbb{R}^3$ as a function of its length parameter $s$ and
\beq
x_1=x(s_0-\epsilon),~~~x_2=x(s_0),~~~x_3=x(s_0+\epsilon)
\eeq
three points on the curve, which at the same time are consecutive corners of an inscribed polygon.   Then the cosine of the angle $\bar{\alpha}$ between the straight edges  at $x_2$ is given by 
\beq
\cos \,\bar{\alpha}~=~\frac{(x_1-x_2)(x_2-x_3)}{\vert x_1-x_2\vert \ \vert x_2-x_3\vert}~=~1-\frac{\kappa^2}{2} \ \epsilon^2+\big (\frac{1}{24}\kappa^4+\frac{1}{12}\kappa^2\tau^2-\frac{1}{12}\kappa \kappa''\big )\ \epsilon ^4+{\cal O}(\epsilon^6)~.\label{metric-alpha}
\eeq
In this expansion use has been made of the  Frenet-Serret formula for the derivatives of the unit vectors $t,n,b$ (note $t=x'(s_0)$) spanning the orthogonal Frenet-Serret frame at $x(s_0)$
\bea
t'&=&\kappa \, n ~,\nonumber\\
n'&=&-\kappa \, t +\tau \, b ~,\nonumber\\
b'&=&-\tau \, n~. \label{Frenet}
\eea
In \eqref{metric-alpha} the leading order $\epsilon^2$ reproduces the curvature, and the next order $\epsilon^4$ is a certain function of curvature and torsion.

There is still another option to get the torsion, and this time directly from the leading term of an expansion. Three consecutive corners of an inscribed
polygon define a plane. Let us then ask for the angle $\bar{\gamma}$ between the normals between two such consecutive planes. To answer this, we look at five consecutive
points $x_1,\dots ,x_5$ on the curve with
\beq
x_1=x(s_0-2\epsilon),~~~x_2=x(s_0-\epsilon),~~~x_3=x(s_0),~~~x_4=x(s_0+\epsilon),~~~x_5=x(s_0+2\epsilon) \label{eps}
\eeq
and find
\beq
\cos \,\bar{\gamma}~=~\frac{(x_2-x_3)\times (x_1-x_3)}{\vert x_2-x_3\vert \ \vert x_1-x_3\vert}\cdot \frac{(x_5-x_3)\times (x_4-x_3)}{\vert x_5-x_3\vert \ \vert x_4-x_3\vert}~=~
1-2\tau^2\ \epsilon^2+{\cal O}(\epsilon^4)~.\label{metric-gamma}
\eeq
\section{The conformal case}
To replace the approximation via ordinary polygons by an approximation which is  conformally covariant, we now describe in some more detail the use
of polygons with circular edges, sketched already in the introduction section. We choose a set of points $\{x_j\}$ on the curve, ordered with respect to increasing
length parameter along the curve. Each triple $(x_{j-1},x_j,x_{j+1})$ fixes a circle $cc_j$. Let us call by $ca^-_j$ its arc from $x_{j-1}$ to $x_j$ and by $ca^+_j$  its arc from $x_j$ to $x_{j+1}$. Then the sequence $( ca^+_1,ca^+_2,ca^+_3,\dots )$ defines an inscribed polygon with circular edges. The same is true for the analogous sequence using the arcs $ca^-_j$, compare fig.1. 

For arbitrary choices of the ordered set of points $\{x_j\}$, under a conformal transformation the inscribed circular polygon is mapped
to the circular polygon inscribed to the image of the original curve and fixed by the image of the original points $\{x_j\}$. However, to make contact
with the conformal invariants for smooth curves \eqref{length-omega}-\eqref{T}, we should consider a situation were consecutive points move closer and closer.
Although metrical length is not conformally invariant, let us start with a setup where consecutive point on the curve have a uniform metrical distance $\epsilon$.

With the choice \eqref{eps} and
\beq
x(s_0+\epsilon)~=~x_3+y_1\, \epsilon+\frac{y_2}{2}\, \epsilon^2+\frac{y_3}{3!}\, \epsilon^3+\dots\label{expan-x}
\eeq
we want to calculate the conformally invariant angle between $ca^+_3$ and $ca^+_2$ at $x_3$, i.e. the cusp angle of the black circular polygon
at the point $x_3$ in fig.1. Using the notation introduced in the appendix, it is given by
\beq
\cos\, \alpha~=~t^{123}_3\cdot t^{234}_3~
\eeq
and via \eqref{tangents} by
\beq
\cos\, \alpha~=~\frac{x_{13}^2\, x_{24}^2+x_{34}^2\, x_{12}^2-x_{23}^2\, x_{14}^2}{2\, x_{12}\, x_{24}\, x_{13}\,x_{34}}~.
\eeq
Then with \eqref{eps} we get, using the Frenet-Serret formula \eqref{Frenet} to express the coefficient vectors $y_j=\frac{d^jx(s)}{ds^j}\vert_{s=s_0}$ in terms of the basis $(t,n,b)$
\beq
\cos\, \alpha~=~1-\frac{1}{8}\, (\kappa^2\tau^2+\kappa'^2)\,  \epsilon^4~+~{\cal O} (\epsilon^5)~.
\eeq
As stressed above, the metrical length is no conformal invariant. Therefore the generic terms in  $\epsilon$-expansions cannot be expected to be  conformal invariants. But there is for sure one exception: the leading nontrivial terms, e.g. the ${\cal O} (\epsilon^4) $ above. Hence
\beq
(\kappa^2\tau^2+\kappa'^2)\, \epsilon^4  
\eeq
is invariant, and we rediscover the differential of the conformal length \eqref{omega}.\\

To find more conformal invariants from higher orders, we now put the points $x_1$ to $x_5$ uniformly equidistant with respect to the conformal length (instead w.r.t. the metrical length). Then all orders in an expansion w.r.t.
this conformal distance will be conformal invariants.

From
\beq
\omega(\epsilon)~=~\int_{s_0}^{s_0+\epsilon}f(s)\, ds,~~~~~f(s)=\sqrt{\nu(s)}=(\kappa^2\tau^2+\kappa'^2)^{\frac{1}{4}}
\eeq
and its related expansion 
\beq
\omega(\epsilon)~=~f(s_0)\,\epsilon+\frac{1}{2}f'\, \epsilon^2+\frac{1}{3!}f''\, \epsilon^3+\dots
\eeq
we get for the inverse function
\beq
\epsilon(\omega)~=~g_1\, \omega+\frac{1}{2} \, g_2\, \omega^2+\frac{1}{3!}\, g_3\, \omega^3+\dots ~,\label{expan-om}
\eeq
with
\beq
g_1=\frac{1}{f}~,~~~~g_j=\frac{g'_{j-1}}{f}~,~~~~j\geq 2~.
\eeq
Now we repeat the calculation of the cusp angle $\alpha$ as an expansion in $\omega$ for the setup
\bea
x_1&=&x(s_0+\epsilon(-2\omega)),~~~x_2~=~x(s_0+\epsilon(-\omega)),~~~x_3~=~x(s_0), \nonumber\\x_4&=&x(s_0+\epsilon(\omega)),~~~x_5~=~x(s_0+\epsilon(2\omega)) \label{omega}~.
\eea
To get the result correctly up to order $\omega^6$ we used in the corresponding Mathematica code the expansions \eqref{expan-x} and \eqref{expan-om}
up to the seventh order and obtained
\beq 
\cos\, \alpha~=~1-\frac{1}{8}\ \omega^4~+~ P\,  \omega^6~+~{\cal O} (\omega^7)~,\label{alpha}
\eeq
with
\bea
P&=&\Big (2 \kappa \tau \kappa'^2 (20 \kappa'\tau' - 
    19\tau\kappa'') +2\kappa^3\tau^3(5\kappa'\tau'-4\tau \kappa'')+\kappa'^2(28\tau^2\kappa'^2-5\kappa''^2+4\kappa'\kappa''')\nonumber\\
    &&~+\kappa^2\big (19\tau^4\kappa'^2-4\kappa'^4 +10\kappa'^2\tau'^2+2\tau\kappa'(2\kappa'\tau''-15\tau'\kappa'')+2\tau^2(5\kappa''^2+2\kappa'\kappa''')\big )\nonumber\\
    &&~+\kappa^4\tau^2(6\tau^4-8\kappa'^2-5\tau'^2+4\tau\tau'')-4\kappa^6\tau^4 \Big)\, \Big (192\,(\kappa^2\tau^2+\kappa'^2)^{\frac{5}{2}}\Big )^{-1}\label{P}~.
    \eea
At first glance this conformal invariant $P$ is not related in an obvious way to conformal curvature $Q$ and torsion $T$. There is an analogous situation to 
the metric case   \eqref{metric-alpha}, where the second nontrivial order is also a certain  function of curvature and torsion.

For the moment, we leave this issue open and turn to the analogue of \eqref{metric-gamma}.  There we considered the angle between the normal vectors
of two consecutive planes,  fixed in both cases by three consecutive points $x_j$. Now we look at the angle between the normal vectors of two
consecutive spheres, fixed both by four consecutive points. Let us denote by $cs(x_1,x_2,x_3,x_4)$ the sphere fixed by the points $x_1,\dots x_4$. Then
the angle between vectors normal to $cs(x_1,x_2,x_3,x_4)$ and $cs(x_2,x_3,x_4,x_5)$ at the point $x_4$ is given by
\beq
\cos\, \gamma~=~\frac{(A_1\, z_1+A_2\, z_2+A_3\, z_3)\cdot (B_2\, z_2+B_3\, z_3+B_5\, z_5)}{\vert A_1\, z_1+A_2\, z_2+A_3\, z_3\vert \, \vert B_2\, z_2
+B_3\, z_3+B_5\, z_5\vert }~.
\eeq
We used the appendix \eqref{sphere-A} and
\beq
B_2=A_2\vert_{x_1\rightarrow x_5}~,~~B_3=A_3\vert_{x_1\rightarrow x_5}~,~~B_5=A_1\vert_{x_1\rightarrow x_5}~.
\eeq
Since our interest concerns the first nontrivial order only, the use of the expansion for $\epsilon (\omega)$ as in \eqref{omega} is not necessary.
Then with the setting
\beq
x_1=x(s_0-3\epsilon),~~~x_2=x(s_0-2\epsilon),~~~x_3=x(s_0-\epsilon),~~~x_4=x(s_0),~~~x_5=x(s_0+\epsilon) \label{sphere-eps}~
\eeq
we obtain
\beq
\cos\, \gamma~=~1-\frac{1}{2}\, T^2\, (\kappa^2\tau^2+\kappa'^2)\,\epsilon^2~+~{\cal O} (\epsilon^3)~=~1-\frac{1}{2}\, T^2\, \omega^2~+~{\cal O} (\omega^3)~. \label{gamma}
\eeq
After this straightforward recovering of the conformal torsion $T$, we come back to the issue of relating the conformal invariant  $P$ from \eqref{P} to $T$ and $Q$.
The denominator in \eqref{P} agrees with the square of that in \eqref{T}, up to the numerical factor $192$. In the nominator there is   $6 \kappa^4\tau^6$, as first term in the third line, which shows up also in the square of the numerator of $T$, again up to a numerical factor. Hence we try to relate
$P-\frac{T^2}{32}$ directly to $Q$ and find \footnote{ A linear combination of $Q$ and $T^2$ also plays a role in some different approach \cite{Langevin}.}
 \beq
 P~=~\frac{1}{24}\, \big (Q+\frac{3}{4}\, T^2\big )~. \label{PQT}
 \eeq
 \\
 
 At the end of this section we still calculate a torsion angle $\beta$ as defined in our papers \cite{Dorn:2020meb,Dorn:2020vzj,Dorn:2023qbj}, and now applied
 to the inscribed polygons under discussion.  Take e.g. for fig.1 the angle between the circumcircle $cc_4$ and the circular edge $ca^+_2$ at the point $x_3$.  
 In fig.1 it is at $x_3$ the angle between the black arc  $ca^+_2$ and the red arc $ca^-_4$
 \beq
 \cos\, \beta~=~t^{123}_3\cdot t^{345}_3~.
\eeq
With \eqref{tangents} from the appendix this yields
\beq
\cos\, \beta~=~\frac{x_{13}^2\, x_{24}^2\,x_{35}^2+ x_{15}^2\, x_{23}^2\, x_{34}^2-x_{13}^2\, x_{25}^2\, x_{34}^2-x_{14}^2\, x_{23}^2\,x_{35}^2}{2\, x_{12}\, x_{13}\,x_{23}\,x_{34}\,x_{35}\,x_{45}}~.
\eeq
For the setting \eqref{omega} this implies
\beq 
\cos\, \beta~=~1-\frac{1}{8}\, T^2\, \omega^6~+~{\cal O} (\omega^7)~.\label{beta}
\eeq

Obviously, we have two options to get $T^2$ from the leading nontrivial term in an $\epsilon$-expansion.  One using spheres for the angle $\gamma$ in \eqref{gamma} and one using circles for the angle $\beta$ in \eqref{beta}. 

But as a crucial message let us stress, that for extracting all three conformal parameters $\omega,Q,T$, one can restrict to calculations of the anglees $\alpha$ and $\beta$  with circles alone \eqref{alpha},\eqref{PQT},\eqref{beta}.
 
\section{Cross ratios and circumcircles}
As a byproduct of our manipulations with angles between crossing circles, we add some comments on generic polygons  \footnote{The inscribed polygons used
in the previous section have edges, which are parts of circles fixed by three consecutive corners. In the generic case the circular edges are not
parts of these circles.}
as considered in our Wilson loop papers  \cite{Dorn:2020meb,Dorn:2020vzj}. For the overall conformal parameterisation of these polygons one needs also cross ratios of the corner points, beyond the 
cusp and torsion angles.  We use the elementary geometry from the appendix to relate cross ratios to angles between certain circumcircles.

We  consider four generic points $x_1,x_2,x_3,x_4 \in  \mathbb{R}^3$. Their conformal properties are characterised by
the cross ratios
\beq
C(i,j,k,l)~=~\frac{ x_{ik}}{ x_{il}}~ \frac{ x_{jl}}{x_{jk}}~.
\eeq
Among them only two are independent. We follow the convention in \cite{Dolan:2000ut} and choose as independent parameters
\beq
u~=~C(1,4,2,3)~=~\frac{ x_{12}}{x_{13} }~ \frac{x_{43}}{ x_{42}}~,~~~v~=~C(1,2,4, 3)~=~\frac{x_{14}}{ x_{13}}~ \frac{ x_{23}}{x_{24}}~.
\eeq

If we take out one of the four points, the remaining three determine a circle, the circumcircle of the related triangle.  In this manner we get four circles $cc^{(i)}$, numerated just by the point omitted.
By definition, each pair of such circles has two points in common, always two points from the set $(x_1,x_2,x_3,x_4)$.  

Let us denote by $t^{(i)}_j$ the unit tangent vector to $cc^{(i]}$ at the point $x_j$. Then the cosine of the crossing angle between e.g. $cc^{(1)}$ and $cc^{(2)}$is given by
\beq
(cc^{(1)},cc^{(2)})~=~t^{(1)}_3t^{(2)}_3~=t^{(1)}_4t^{(2)}_4~.
\eeq

Our aim is now to express all the six $(cc^{(i)},cc^{(j)})$ in terms of the two independent cross ratios $u$ and $v$.
Then with \eqref{tangents} and suitable adaption of indices we get 
\bea
p=\cos\,\varphi&:=&(cc^{(1)},cc^{(2)})~=~(cc^{(3)},cc^{(4)})~=~\frac{1+v^2-u^2}{2\ v}~,\nonumber\\
q=\cos\,\psi&:=&(cc^{(1)},cc^{(4)})~=~(cc^{(2)},cc^{(3)})~=~\frac{1+u^2-v^2}{2\ u}~,\nonumber\\
r=\cos\,\chi&:=&(cc^{(1)},cc^{(3)})~=~(cc^{(2)},cc^{(4)})~=~\frac{1-u^2-v^2}{2\ u\  v}~.\label{abc,uv}
\eea
\newpage
\noindent
Using the above parametrisation of $p,q$ and $r$ by $u$ and $v$ we obtain
\beq 
r~=~p\, q-\sqrt{(1-p^2)(1-q^2)}~,~\mbox{i.e.}~~~ \cos\,\chi=\cos\,(\varphi +\psi )~.\label{r,pq}
\eeq
\begin{figure}[h!]
\begin{center}
  \includegraphics[width=9cm]{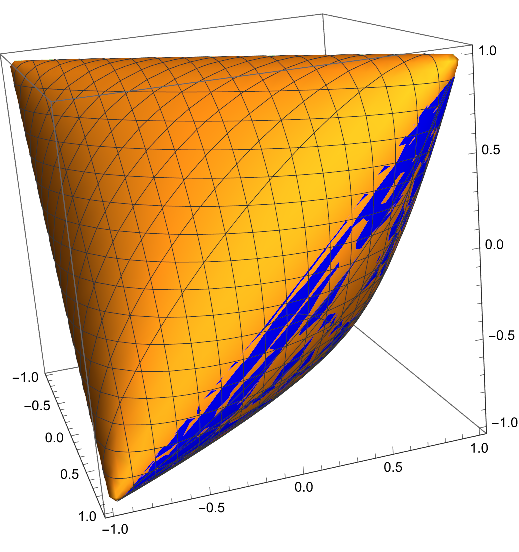}
  \end{center}

\caption {\it In yellow the boundary of the rounded tetrahedron \eqref{sym-constraint} and in yellow/blue that part, for which  $(p,q,r)$  is allowed by  \eqref{abc,uv} and  
\eqref{uv-constraint}.}
\label{tetra}
\end{figure}
A manifest symmetric form of the constraint\footnote{Interestingly, this cubic  constraint for three parameters appears also in the context of the conformal geometry of triangles with circular edges  \cite{Dorn:2020vzj} and the conformal properties of null hexagons \cite{Dorn:2012cn},\cite{Alday:2009dv}.} between the three angles is
\beq
p^2+q^2+r^2-2p\, q\, r~=~1~.\label{sym-constraint}
\eeq
However, one has to keep in mind, that not all points  $(p,q,r)$  allowed by \eqref{sym-constraint} belong to a possible geometric setup for $x_1,x_2,x_3,x_4 \in  \mathbb{R}^3$.  First of all,  \eqref{sym-constraint} includes also points corresponding to the positive sign in front of the square root in \eqref{r,pq}. Furthermore, since $p,q,r$ as cosines of angles have to have a modulus smaller or equal to 1, the equations \eqref{abc,uv} imply for the by definition
positive cross ratios $u,v$
\bea
u\geq 0~,&&v\geq 0~,\nonumber\\
u+v\geq1~,&&\vert u-v\vert\leq 1~.\label{uv-constraint}
\eea
Of course these constraints on the cross ratios are well-known. A straightforward way to derive them without consideration of our circles is given
by conformal mapping e.g. $x_4$ to infinity. Then  \eqref{uv-constraint} is a consequence of the triangle 
inequalities for the triangle with corners $x_1,x_2,x_3$.

The surface in $(p,q,r)$-space, defined by \eqref{sym-constraint} and $\vert p\vert,\vert q\vert,\vert r\vert \leq 1$ looks like the boundary of a tetrahedron with rounded edges. Only one of the four faces is allowed by \eqref{abc,uv} and  \eqref{uv-constraint}. It is just the face opposite to the vertex $(p,q,r)=(-1,-1,1)$, see fig.\ref{tetra}. 

The vertices $(1,1,1),~(1,-1,-1)$ and $(-1,1,-1)$ are via \eqref{abc,uv} the images of the boundary of the allowed region \eqref{uv-constraint} in the $(u,v)$-plane: $u+v=1,~u-v=1$ and $v-u=1$, respectively.  They correspond to the limiting case, where all four points are located
on just one circle or one straight line. All crossing angles are then equal to $0$ or $\pi$,  depending on ordering issues generated by the sign convention in \eqref{tangents}.
\\[20mm]
{\bf Acknowledgement}\\[2mm]
I thank the Quantum Field and String Theory Group at Humboldt University for kind hospitality.
\section*{Appendix}
Three points $x_1,x_2,x_3$ fix a circle uniquely. Its unit tangent vectors at these points are given by
\bea
t^{(123)}_1&=& \frac{1}{x_{23}}\big ( \frac{x_{13}}{x_{12}}(x_2-x_1)+\frac{x_{12}}{x_{13}}(x_1-x_3)\big )~,  \nonumber\\
t^{(123)}_2&=&  \frac{1}{x_{13}}\big ( \frac{x_{12}}{x_{23}}(x_3-x_2)+\frac{x_{23}}{x_{12}}(x_2-x_1)\big )~,  \nonumber\\
t^{(123)}_3&=& \frac{1}{x_{12}}\big ( \frac{x_{23}}{x_{13}}(x_1-x_3)+\frac{x_{13}}{x_{23}}(x_3-x_2)\big )~,\label{tangents}
\eea 
with $x_{ij}=\vert x_i-x_j\vert$.

Four points  $x_1,x_2,x_3,x_4$ fix a sphere uniquely. Its center $x_c$ is given by
\beq
x_c=x_4 +\frac{1}{D}\, \big ( A_1\, (x_1-x_4)+A_2\, (x_2-x_4)+A_3\, (x_3-x_4)\big )~,\label{sphere}
\eeq
with $(z_j=x_j-x_4)$
\bea
A_1&=&\left \vert  \begin{array}{c} z_1z_1,\, z_1z_2,\, z_1z_3\\ z_2 z_2,\, z_2z_2,\, z_2z_3\\ z_3z_3,\, z_3z_2,\, z_3z_3\end{array}\right \vert~,~~~A_2~=~\left \vert  \begin{array}{c} z_1z_1,\, z_1z_1,\, z_1z_3\\ z_2 z_1,\, z_2z_2,\, z_2z_3\\ z_3z_1,\, z_3z_3,\, z_3z_3\end{array}\right \vert ~,\label{sphere-A}\\[2mm]
A_3&=&\left \vert  \begin{array}{c} z_1z_1,\, z_1z_2,\, z_1z_1\\ z_2 z_1,\, z_2z_2,\, z_2z_2\\ z_3z_1,\, z_3z_2,\, z_3z_3\end{array}\right \vert~,~~~~D~=~\left \vert  \begin{array}{c} z_1z_1,\, z_1z_2,\, z_1z_3\\ z_2 z_1,\, z_2z_2,\, z_2z_3\\ z_3z_1,\, z_3z_2,\, z_3z_3\end{array}\right \vert~. \nonumber
\eea

\end{document}